\documentclass{article}
\usepackage[margin=1in]{geometry} 
\usepackage{amsmath, amsthm, amssymb, amsfonts, fancyhdr, comment, parskip, graphicx, xcolor}
\usepackage[colorlinks]{hyperref}
\usepackage{mathtools}
\usepackage{upgreek}
\usepackage{subcaption}
\usepackage{bm}
\usepackage{tikz, lmodern}
\usepackage{enumerate}

\usetikzlibrary{graphs}
\usetikzlibrary{shadings}
\usetikzlibrary{shapes.geometric}
\usetikzlibrary{arrows}

\tikzstyle{vertex} = [fill,shape=circle,node distance=80pt]
\tikzstyle{edge} = [fill,opacity=0.5,fill opacity=0.5,line cap=round, line join=round, line width=50pt]
\tikzstyle{mededge} = [fill,opacity=0.5,fill opacity=0.5,line cap=round, line join=round, line width=57.5pt]
\tikzstyle{thickedge} = [fill,opacity=0.5,fill opacity=0.5,line cap=round, line join=round, line width=65pt]
\tikzstyle{elabel} =  [fill,shape=circle,node distance=30pt]

\pgfdeclarelayer{background}
\pgfsetlayers{background,main}

\newtheorem{theorem}{Theorem}[section]
\newtheorem{corollary}[theorem]{Corollary}
\newtheorem{lemma}[theorem]{Lemma}
\newtheorem {conjecture}{Conjecture}
\theoremstyle{remark}
\newtheorem{observation}{Observation}
\theoremstyle{definition}
\newtheorem{definition}[theorem]{Definition}

\newcommand{\tu}[2]{\tau^{(#1)}(#2)}
\newcommand{\vu}[2]{\nu^{(#1)}(#2)}

\title{New bounds on a generalization of Tuza's conjecture}
\author{Alex Parker \thanks{Department of Mathematics, Iowa State University, Ames, IA.\,\,\,\,Email: \texttt{abparker@iastate.edu}}}
\date{\today}

\begin{document}

\maketitle

\begin{abstract}
    For a $k$-uniform hypergraph $H$, let $\nu^{(m)}(H)$ denote the maximum size of a set $S$ of edges of $H$ whose pairwise intersection has size less than $m$. Let $\tau^{(m)}(H)$ denote the minimum size of a set $S$ of $m$-sets of $V(H)$ such that every edge of $H$ contains some $m$-set from $S$. A conjecture by Aharoni and Zerbib, which generalizes a conjecture of Tuza on the size of minimum edge covers of triangles of a graph, states that for a $k$-uniform hypergraph $H$, $\tau^{(k - 1)}(H)/\nu^{(k - 1)}(H) \leq \left \lceil \frac{k + 1}{2} \right \rceil$. In this paper, we show that this generalization of Tuza's conjecture holds when $\nu^{(k - 1)}(H) \leq 3$. As a corollary, we obtain a graph class which satisfies Tuza's conjecture. We also prove various bounds on $\tau^{(m)}(H)/\nu^{(m)}(H)$ for other values of $m$ as well as some bounds on the fractional analogues of these numbers.
\end{abstract}

\section{Introduction}
\subsection{Definitions and Notation}
Throughout this paper, unless otherwise specified, we will only be concerned with $k$-uniform hypergraphs for $k \geq 3$. We start by establishing some definitions and notation which will be used throughout the paper.

For a set $S$ with $x \in S$, $y \not \in S$, we denote $S \setminus \{x\}$ by $S - x$ and $S \cup \{y\}$ by $S + y$. For a set $Z$ with $|Z| = 2$, when we say $z \in Z$, we will let $\overline{z} = Z - z$. For a hypergraph $H$, we will use both $E(H)$ and $H$ to mean the edge set of $H$.
Let $H$ be a $k$-uniform hypergraph with vertex set $V$ and edge set $E$. A \emph{matching} of $H$ is any collection of disjoint edges of $H$. We denote the largest matching of $H$ by $\nu(H)$. A \emph{cover} of $H$ is a set $C \subseteq V$ such that for every $e \in E$, there is some $v \in C \cap e$. We denote the size of the smallest cover of $H$ by $\tau(H)$. Clearly, for any $k$-uniform hypergraph $H$, $\nu(H) \leq \tau(H) \leq k\nu(H)$. \\

These definitions may be generalized in the following way: for $1 \leq m \leq k - 1$, an \emph{$m$-matching} of $H$ is a collection $M$ of edges of $H$ such that for any $e, e' \in M$, $|e \cap e'| < m$. We denote the size of the largest $m$-matching of $H$ by $\nu^{(m)}(H)$. Observe that $\nu(H) = \nu^{(1)}(H)$. An \emph{$m$-cover} of $H$ is a set $C \subseteq \binom{V}{m}$ such that for every $e \in H$, there is some $c \in C$ with $c \subseteq e$. We denote the size of the smallest $m$-cover of $H$ by $\tau^{(m)}(H)$. Again, observe that $\tau(H) = \tau^{(1)}(H)$. Similar to the inequality above, we trivially have $\nu^{(m)}(H) \leq \tau^{(m)}(H) \leq \binom{k}{m} \nu^{(m)}(H)$. The main aim of this paper will be to improve the ratio $\tau^{(m)}(H)/\nu^{(m)}(H)$ for various values of $m$ and $\nu^{(m)}(H)$. \\

We will also study the fractional versions of these parameters. A \emph{fractional $m$-matching} is a function $f: E(H) \to \mathbb{R}_{\geq 0}$ such that for every $S \in \binom{V}{m}$, $\sum_{e \supseteq S} f(e) \leq 1$. The \emph{size} of a fractional $m$-matching is $|f| = \sum_{e \in E(H)} f(e)$. A \emph{fractional $m$-cover} is a function $c: \binom{V}{m} \to \mathbb{R}_{\geq 0}$ such that for every $e \in H$, $\sum_{S \in \binom{e}{m}} c(S) \geq 1$. The \emph{size} of a fractional $m$-cover is $|c| = \sum_{S \in \binom{V}{m}} c(S)$. The fractional $m$-matching number, $\nu^{*(m)}(H)$, and the fractional $m$-cover number $\tau^{*(m)}(H)$ are defined to be the maximum size of a fractional $m$-matching and the minimum size of a fractional $m$-cover, respectively. We will denote $\nu^{*(1)}(H)$ by $\nu^*(H)$ and $\tau^{*(1)}(H)$ by $\tau^*(H)$. Observe that by LP duality, we always have $\nu^{*(m)}(H) = \tau^{*(m)}(H)$. Also, observe that an $m$-matching is a fractional $m$-matching and an $m$-cover is a fractional $m$-cover. For any $k$-uniform hypergraph $H$ and $1 \leq m \leq k - 1$, we have:
\begin{equation*}
    \nu^{(m)}(H) \leq \nu^{*(m)}(H) = \tau^{*(m)}(H) \leq \tau^{(m)}(H) \leq \binom{k}{m} \nu^{(m)}(H).
\end{equation*}

\subsection{A generalization of Tuza's conjecture}
We introduce some notation which will be used throughout the paper. Let $\mathcal{H}_k$ denote the family of all $k$-uniform hypergraphs. Then, define the following functions:
\begin{center}
    \begin{itemize}
        \item $h(k, m) = \sup \left \{ \frac{\tau^{(m)}(H)}{\nu^{(m)}(H)} : H \in \mathcal{H}_k \right \}$ \\
        \item $g_i(k, m) = \sup \left \{ \frac{\tau^{(m)}(H)}{\nu^{(m)}(H)} : H \in \mathcal{H}_k \text{ and } \nu^{(m)}(H) = i \right \}$ \\
        \item $h^*(k, m) = \sup \left \{ \frac{\tau^{*(m)}(H)}{\nu^{(m)}(H)} : H \in \mathcal{H}_k \right \}$ \\
        \item $g_i^*(k, m) = \sup \left \{ \frac{\tau^{*(m)}(H)}{\nu^{(m)}(H)} : H \in \mathcal{H}_k \text{ and } \nu^{(m)}(H) = i \right \}.$
    \end{itemize}
\end{center}

For reference, some previous papers used $g(k, m)$ for $g_1(k,m)$ and $g^*(k, m)$ for $g_1^*(k, m)$. Observe that by definition, we have:
\begin{align*}
    g_i^*(k, m) \leq g_i(k, m) &\leq h(k, m) \\
    g_i^*(k, m) \leq h^*(k, m) &\leq h(k, m).
\end{align*}

A famous conjecture of Tuza~\cite{tuza90} states that for any graph $G$, the minimum number of edges needed to intersect every triangle in $G$ ($\tau_t(G)$) is at most twice the maximum number of edge disjoint triangles in $G$ ($\nu_t(G)$). If true, this conjecture is tight as seen e.g., by $K_4$ or $K_5$. The conjecture has been shown to be true for various families of graphs (see e.g.~\cite{botler19}, \cite{tuza90}). Haxell~\cite{haxell99} proved the best known general upper bound of $\tau_t(G) \leq \frac{66}{23} \nu_t(G)$. \\

Note that for a graph $G$, if we define the triangle graph of $G$, $T(G)$, to be the hypergraph with edges corresponding to the triangles of $G$, Tuza's conjecture states that for any graph $G$, $\tau^{(2)}(T(G))/\nu^{(2)}(T(G)) \leq 2$. A conjecture of Aharoni and Zerbib generalizes Tuza's, conjecturing that for all $3$-uniform hypergraphs $H$, $\tau^{(2)}(H)/\nu^{(2)}(H) \leq 2$ (i.e. $h(3,2) \leq 2$). \\

Furthermore, they conjectured that a similar bound should hold for hypergraphs of any fixed uniformity:
\begin{conjecture}[\cite{az20}]
    Let $k \geq 3$. Then, $h(k, k - 1) \leq \left \lceil \frac{k + 1}{2} \right \rceil$.
\end{conjecture}

Again, if true, the bound is tight as seen by the following example from~\cite{az20}: for $H = \binom{[k + 1]}{k}$, the $k$-uniform hypergraph containing all $k$-subsets of $[k+1]$, one can easily check that $\nu^{(k - 1)}(H) = 1$ and $\tau^{(k - 1)}(H) = \left \lceil \frac{k + 1}{2} \right \rceil$.

\subsection{The paper}
We begin by studying the function $g_i(k, k - 1)$ in section 2. In~\cite{az20}, Aharoni and Zerbib showed that $g_1(k, k - 1) \leq \left \lceil \frac{k + 1}{2} \right \rceil$. We prove the same bound for $g_2(k, k - 1)$ and $g_3(k, k - 1)$:
\begin{theorem}\label{thm:main1}
    Let $H$ be a $k$-uniform hypergraph with $\nu^{(k - 1)}(H) = 2$. Then,
    \begin{equation*}
        \tau^{(k - 1)}(H) \leq 2 \left \lceil \frac{k + 1}{2} \right \rceil.
    \end{equation*}
\end{theorem}

\begin{theorem}\label{thm:main2}
    Let $H$ be a $k$-uniform hypergraph with $\nu^{(k - 1)}(H) = 3$. Then,
    \begin{equation*}
        \tau^{(k - 1)}(H) \leq 3 \left \lceil \frac{k + 1}{2} \right \rceil.
    \end{equation*}
\end{theorem}

This immediately implies the following:
\begin{corollary}
    Let $G$ be a graph with the property that $G$ does not contain $4$ edge-disjoint triangles. Then, Tuza's conjecture holds for $G$.
\end{corollary}

In section 3, we study $g_1(k, m)$ for various values of $m$. We prove the first non-trivial upper bounds for $g_1(k, m)$ when $\frac{k}{2} \leq m \leq k - 2$.
\begin{theorem}\label{thm:mGEQk2}
    Let $k \geq 6$ and suppose $\frac{k}{2} \leq m \leq k - 2$. Then, $g_1(k, m) \leq \binom{k}{m} - m$.
\end{theorem}

\begin{theorem}\label{thm:mEqkMinus2}
    Let $k \geq 3$. Then, we have:
    \begin{equation*}
        g_1(k, k - 2) \leq \left \lceil \frac{k^2}{4} \right \rceil =
            \begin{cases}
                \frac{1}{4}(k^2 + 3), & \text{if } k \text{ odd}, \\
                \frac{1}{4}k^2, & \text{if } k \text{ even}.
            \end{cases}
    \end{equation*}
\end{theorem}

Aharoni and Zerbib~\cite{az20} previously showed that $g_1(k, 2) < \binom{k}{2}$ and $g_1(4, 2) = 4$. We go on to improve the upper bound of $g_1(5, 2)$ (the first remaining open case when $m = 2$) with the best previous bound being $g_1(5, 2) \leq 9$.

\begin{theorem}\label{thm:g52}
    We have $6 \leq g_1(5, 2) \leq 7$.
\end{theorem}

The lower bound has not been mentioned in previous papers but comes from the $2$-cover number of the (unique) symmetric $2-(11,5,2)$ design (an explicit construction can be seen in Table 1.19 in~\cite{colburn}).

In section 4, we study the fractional variants of the problem and prove bounds on $g_1^*(k, m)$ for certain choices of $m$:


\begin{theorem}\label{thm:fracG2kk}
    For all $k \geq 2$, $g_1^*(2k, k) \leq \left( \frac{1}{2} + \frac{1}{2(k + 1)}\right)\binom{2k}{k}$.
\end{theorem}

The proof of this theorem is followed by a lemma, generalizing a result from~\cite{basit22}, that allows us to obtain upper bounds on $h^*(k,m)$ from upper bounds on $g_1^*(k,m)$. When $m = k/2$, this gives the following corollary:

\begin{corollary}\label{cor:fracH2kk}
    For all $k \geq 2$, $h^*(2k, k) \leq \left(1 - \frac{k}{4(k+2)} \right) \binom{2k}{k}$.
\end{corollary}

We also prove a fractional upper bound on $g_1^*(k, k-2)$ from which a bound for $h^*(k, k-2)$ may be derived in the same manner as above.

\begin{theorem}\label{thm:fracGkkMinus2}
    $g_1^*(k, k - 2) \leq \frac{1}{6} \binom{k - 2}{2} + 2k - 3$.
\end{theorem}

It should be noted that other fractional variations and results have been shown in~\cite{basit22}, \cite{gur22}, among others.

\section{$g_i(k, k - 1)$}

We begin this section with some useful definitions and a few short lemmas.

\begin{definition}
    Let $H$ be a $k$-uniform hypergraph and $M$ be a maximum $(k - 1)$-matching in $H$. For any vertex $v \in V(H)$, we denote $d_M(v)$ to be the number of edges of $M$ that contain $v$. For each $e \in M$, define the following two sets:
    \begin{align*}
        S_e &=  \{h \in H : |e \cap h| \geq k - 1 \text{ and } |h \cap f| < k - 1 \text{ for all } f \in M - e\} \\
        T_e &=  \{h \in H : |e \cap h| \geq k - 1\}.
    \end{align*}
\end{definition}

\begin{lemma}\label{lem:disjS_i}
    Let $H$ be a $k$-uniform hypergraph and $M$ a maximum $(k - 1)$-matching in $H$. Then, for any $e, f \in M$, $S_e \cap S_f = \emptyset$. Further, $\nu^{(k - 1)}(S_e) = 1$, which implies $\tau^{(k - 1)}(S_e) \leq g_1(k, k - 1)$.
\end{lemma}

\begin{proof}
    This follows directly from the definition of $S_e$.
\end{proof}

\begin{lemma}\label{lem:T_iLB}
    Let $H$ be a $k$-uniform hypergraph and let $M$ be a maximum $(k - 1)$-matching in $H$. If there exists some $e \in M$ such that $\tau^{(k - 1)}(T_e) \leq \left \lceil \frac{k + 1}{2} \right \rceil$, then
    \begin{equation*}
        \tu{k - 1}{H} \leq \left \lceil \frac{k + 1}{2}\right \rceil + (\vu{k - 1}{H} - 1)g_{\vu{k - 1}{H} - 1}(k, k - 1).
    \end{equation*}
\end{lemma}

\begin{proof}
    Let $H$ be a $k$-uniform hypergraph and let $M$ be a maximum $(k - 1)$-matching in $H$. Suppose there exists some $e \in M$ such that $\tu{k - 1}{T_e} \leq \left \lceil \frac{k + 1}{2} \right \rceil$. We claim that $H - T_e$ has matching number at most $\vu{k - 1}{H} - 1$. Suppose not. Then, there exists some matching $M'$ of $H - T_e$ of size at least $\vu{k - 1}{H}$. By definition, all edges of $H - T_e$ intersect $e$ in at most $k - 2$ vertices. But then, $M' + e$ is a larger matching than $M$, a contradiction. Therefore, we have:
    \begin{equation*}
        \tu{k - 1}{H} \leq \tu{k - 1}{T_e} + \tu{k - 1}{H - T_e} \leq \left \lceil \frac{k + 1}{2} \right \rceil + (\vu{k - 1}{H} - 1)g_{\vu{k - 1}{H} - 1}(k, k - 1).
    \end{equation*}
\end{proof}

\begin{lemma}\label{lem:connected}
    Let $H$ be a $k$-uniform hypergraph and let $M$ be a maximum $(k - 1)$-matching in $H$. If there exists a partition $P_1, P_2$ of the edges of $M$ such that for all $e \in P_1$ and $e' \in P_2$, $|e \cap e'| < k - 2$, then $T_e \cap T_{e'} = \emptyset$ and
    \begin{equation*}
        \tu{k - 1}{H} \leq |P_1| g_{|P_1|}(k, k - 1) + |P_2| g_{|P_2|}(k, k - 1).
    \end{equation*}

    We call such a matching disconnected.
\end{lemma}

\begin{proof}
    Let $H$ be a $k$-uniform hypergraph and let $M$ be a maximum $(k - 1)$-matching in $H$. Suppose there exists a partition $P_1, P_2$ of the edges of $M$ such that for all $e \in P_1$ and $e' \in P_2$, $|e \cap e'| < k - 2$. Now, let $e \in P_1$, $e' \in P_2$ and suppose $f \in T_e$. Then, $f$ intersects $e$ in $k - 1$ vertices and therefore, $f$ can only intersect $e'$ in at most $k - 2$ vertices. So, $T_e \cap T_{e'} = \emptyset$. This means that the edges of $H$ are the disjoint union of the sets $H_1 := \bigcup_{e \in P_1} T_e$ and $H_2 := \bigcup_{e' \in P_2} T_{e'}$. Also, because there is no intersection of size $k - 1$ between any edge of $H_1$ and any edge in $P_2$, $\vu{k - 1}{H_1} = |P_1|$. Similarly, $\vu{k - 1}{H_2} = |P_2|$. Therefore,
    \begin{equation*}
        \tu{k - 1}{H} \leq \tu{k - 1}{H_1} + \tu{k - 1}{H_2} \leq |P_1| g_{|P_1|}(k, k - 1) + |P_2| g_{|P_2|}(k, k - 1).
    \end{equation*}
\end{proof}

\begin{lemma}\label{lem:pendantEdge}
    Let $H$ be a $3$-uniform hypergraph and let $M$ be a maximum $2$-matching in $H$. If there exists some $e \in M$ such that $\sum_{v \in e} d_{M}(v) \leq 4$ and $\tu{2}{S_e} = 1$, then
    \begin{equation*}
        \tu{2}{H} \leq 4 + (\vu{2}{H} - 2)g_{\vu{2}{H} - 2}(k, k - 1).
    \end{equation*}
\end{lemma}

\begin{proof}
    Let $H$ be a $3$-uniform hypergraph and let $M$ be a maximum $2$-matching in $H$. Suppose there exists some $e \in M$ such that $\sum_{v \in e} d_{M}(v) \leq 4$ and $\tu{2}{S_e} = 1$. This means that there are two vertices in $e$ not contained in any other edge of $M$ and at most one vertex of $e$ contained in at most one other edge, say $f$, of $M$. Then, it is clear that $(T_e - S_e) \subseteq T_f$. Furthermore, $\vu{2}{H - T_e - T_f} = |M| - 2$. Otherwise, if we may find a $2$-matching $M'$ of $H - T_e - T_f$ of size greater than $|M| - 2$, then $M' + e + f$ is larger than $M$, a contradiction. Now, let $S$ be a $2$-set, which $2$-covers $S_e$. Then, since $T_e - S_e \subseteq T_f$ and $S_f \subseteq T_f$, taking $\binom{f}{2}$ to $2$-cover $T_f$, we have found a $2$-cover of $T_e \cup T_f$ of size $4$. Therefore, we have:
    \begin{equation*}
        \tu{2}{H} \leq \tu{2}{T_e \cup T_f} + \tu{2}{H - T_e - T_f} \leq 4 + (\vu{2}{H} - 2)g_{\vu{2}{H} - 2}(k, k - 1).
    \end{equation*}
\end{proof}

We now refine the $\nu^{(k - 1)} = 1$ result of Aharoni and Zerbib~\cite{az20} in order to help with our proof of the $\nu^{(k - 1)} \in \{2, 3\}$ cases. First, we reiterate a lemma from~\cite{basit22}:

\begin{lemma}[Lemma 2.2 from~\cite{basit22}] \label{lem:uniqVtx}
    Let $H$ be a $k$-uniform hypergraph with $\vu{k - 1}{H} = 1$. Then, either $\tu{k-1}{H} = 1$ or for any edge $e \in E(H)$, there exists a unique vertex $v \in V(H) - V(e)$ such that for all $e' \in E(H) - e$, $e' - e = \{v\}$.
\end{lemma}

Now, we are ready to refine the $\nu^{(k - 1)} = 1$ result from~\cite{az20}.

\begin{lemma}\label{lem:trueTau}
    Let $H$ be a $k$-uniform hypergraph with $\vu{k-1}{H} = 1$. Then, either $\tu{k-1}{H} = 1$ or $\tu{k-1}{H} \leq \left \lceil \frac{e(H)}{2} \right \rceil$.
\end{lemma}

\begin{proof}
    Let $H$ be a $k$-uniform hypergraph with $k \geq 3$. Suppose $\vu{k-1}{H} = 1$ and $\tu{k-1}{H} \neq 1$. Let $e \in E(H)$ and let $v \in V(H) - V(e)$ be the vertex described in Lemma~\ref{lem:uniqVtx}. Let $e_1, \dots, e_{e(H) - 1}$ denote the edges of $H - e$. Observe that for any $1 \leq i \neq j \leq e(H) - 1$, $|e_i \cap e_j \cap e| = k - 2$. \\
    Suppose $e(H)$ is odd. For $1 \leq i \leq \frac{e(H) - 1}{2}$, we may cover $e_{2i - 1}, e_{2i}$ with the set $(e_{2i - 1} \cap e_{2i} \cap e) + v$. Then, we may cover $e$ with any set from $\binom{e}{k - 1}$, giving a $(k - 1)$-cover of size $\frac{e(H) - 1}{2} + 1 = \frac{e(H) + 1}{2} = \left \lceil \frac{e(H)}{2} \right \rceil$. \\
    Suppose $e(H)$ is even. For $1 \leq i \leq \frac{e(H) - 2}{2}$, we may cover $e_{2i - 1}, e_{2i}$ with the set $(e_{2i - 1} \cap e_{2i} \cap e) + v$. Then, we may cover $e_{e(H) - 1}, e$ with the set $e_{2i - 1} \cap e$, giving a $(k - 1)$-cover of size $\frac{e(H) - 2}{2} + 1 = \frac{e(H)}{2} = \left \lceil \frac{e(H)}{2} \right \rceil$.
\end{proof}

We obtain the $\nu^{(k - 1)} = 1$ result as a corollary:
\begin{corollary}\label{cor:nuEquals1}
    We have $g_1(k, k - 1) \leq \left \lceil \frac{k + 1}{2} \right \rceil$.
\end{corollary}

\begin{proof}
    Let $H$ be a $k$-uniform hypergraph with $\vu{k-1}{H} = 1$. We may assume $\tu{k-1}{H} > 1$. Let $e \in H$ and let $v \in V(H) - V(e)$ be the unique vertex as described in Lemma~\ref{lem:uniqVtx}. Now, aside from $e$, every other edge of $H$ consists of $v$ together with some $(k-1)$-subset of $e$. Since $e$ has $k$ different $(k - 1)$-subsets, then the total number of edges of $H$ is at most $k + 1$. The result now follows from Lemma~\ref{lem:trueTau}.
\end{proof}

Next, we prove the case when $\nu^{(k-1)} = 2$.

\begin{proof}[Proof of Theorem~\ref{thm:main1}]
    Let $H$ be a $k$-uniform hypergraph with $k \geq 3$. Suppose $\vu{k-1}{H} = 2$. If there exists a $(k - 1)$-matching of $H$, $\{e, f\}$, where $|e \cap f| < k - 2$, then $e, f$ is a disconnected matching and we are done by Lemma~\ref{lem:connected} together with Lemma~\ref{lem:trueTau}.

    Suppose then that for any maximum $(k - 1)$-matching $\{e, f\}$ in $H$, $|e \cap f| = k - 2$. To this end, let $\{e, f\}$ be a $(k - 1)$-matching of $H$ with
    \begin{align*}
        e &=    S \cup \{u_1, u_2\} \\
        f &=    S \cup \{v_1, v_2\}.
    \end{align*}

    Here, $S = e \cap f$ is a $(k - 2)$-subset of $V(H)$. Since $\vu{k-1}{S_e} = \vu{k-1}{S_f} = 1$, then as noted before, $\tu{k-1}{S_e} \leq \left \lceil \frac{k + 1}{2} \right \rceil$ and $\tu{k-1}{S_f} \leq \left \lceil \frac{k + 1}{2} \right \rceil$. If every edge of $S_e$ contains $S$, then we may cover $T_e$ with the sets $S + u_1$ and $S + u_2$. Next, we may cover $S_f$ with at most $\left \lceil \frac{k + 1}{2} \right \rceil$ $(k - 1)$-sets, giving a cover of $H$ of size at most
    \begin{equation*}
        2 + \left \lceil \frac{k + 1}{2} \right \rceil \leq 2 \left \lceil \frac{k + 1}{2} \right \rceil.
    \end{equation*}
    
    Similarly, we may find a cover of suitable size if every edge of $S_f$ contains $S$. Further, if $\tu{k-1}{S_e} = 1$, then we may cover $S_e$ with one $(k - 1)$-set and cover the rest of $H$ with elements from $\binom{f}{k - 1}$, giving a cover of size at most
    \begin{equation*}
        1 + k \leq 2 \left \lceil \frac{k + 1}{2} \right \rceil.
    \end{equation*}
    
    Similarly, we may find a cover of suitable size if $\tu{k-1}{S_f} = 1$. So, we may assume $\tu{k-1}{S_e} \neq 1$, $\tu{k-1}{S_f} \neq 1$ and that there exists $e' \in S_e - e$, $f' \in S_f - f$ such that $S \not \subseteq e'$ and $S \not \subseteq f'$. \\
    
    If the unique vertex for all edges of $S_e - e$ described in Lemma~\ref{lem:uniqVtx} is not contained in $f - e$, then $e', f$ is a disconnected matching. So, we may assume that for all $e'' \in S_e - e$, $e'' - e = v$ for some $v \in \{v_1, v_2\}$. By a symmetric argument, for all $f'' \in S_f - f$, $f'' - f = u$ for some $u \in \{u_1, u_2\}$. \\

    This tells us that every edge in $S_e - e$ is of the form $S' \cup \{u_1, u_2, v\}$ for some $S' \in \binom{S}{k - 3}$ (i.e. there are at most $k - 1$ edges in $S_e$). Similarly, every edge in $S_f - f$ is of the form $S'' \cup \{v_1, v_2, u\}$ for some $S'' \in \binom{S}{k - 3}$ (i.e. there are at most $k - 1$ edges in $S_f$). By Lemma~\ref{lem:trueTau}, we may cover every edge in $S_e$ with at most $\left \lceil \frac{k - 1}{2} \right \rceil$ $(k - 1)$-sets and we may cover every edge in $S_f - f$ with at most $\left \lceil \frac{k - 2}{2} \right \rceil$ $(k - 1)$-sets. Finally, we may cover $T_f - S_f + f$ with the sets $S + v_1$ and $S + v_2$. This gives us a cover of $H$ of size at most
    \begin{equation*}
        \left \lceil \frac{k - 1}{2} \right \rceil + \left \lceil \frac{k - 2}{2} \right \rceil + 2 = k + 1 \leq 2 \left \lceil \frac{k + 1}{2} \right \rceil.
    \end{equation*}
\end{proof}

Now, we are ready to prove the $\nu^{(k- 1)} = 3$ case:

\begin{proof}[Proof of Theorem~\ref{thm:main2}]
    We break the proof into two parts. In the first part, we assume we are dealing with a $3$-uniform hypergraph. In the second part, we will deal with an arbitrary $k$-uniform hypergraph with $k \geq 4$. \\
    Let $H$ be a $3$-uniform hypergraph and let $M = \{e, f, g\}$ be a maximum $2$-matching in $H$. If $M$ is disconnected, then the result follows from Lemma~\ref{lem:connected}. So, suppose $M$ is connected. We may assume $|e \cap f| = 1$ and $|e \cap g| = 1$. Then, $M$ looks like one of the matchings from Figure~\ref{fig:3UM}. \\

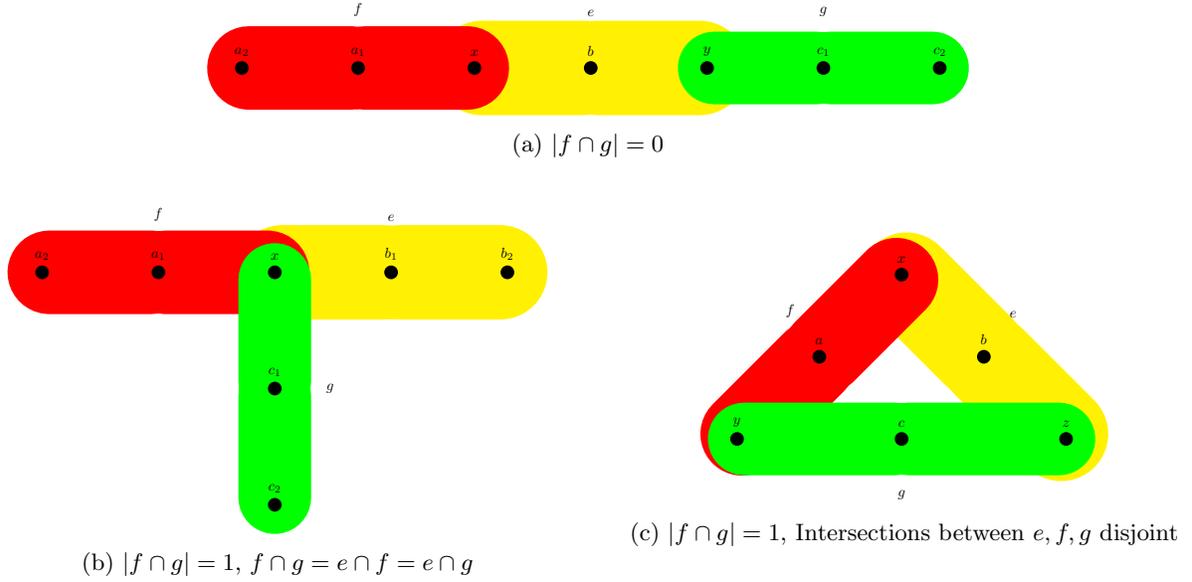
\begin{figure}[h]
    \begin{subfigure}[h]{\textwidth}
        \centering
        \scalebox{0.55}
        {
            \begin{tikzpicture}
                [every label/.append style={font=\LARGE}]
                \node[vertex, label=above:\(a_2\)] (a2) {};
                \node[vertex,right of=a2,label=above:\(a_1\)] (a1) {};
                \node[vertex,right of=a1,label=above:\(x\)] (x) {};
                \node[vertex,right of=x,label=above:\(b\)] (b) {};
                \node[vertex,right of=b,label=above:\(y\)] (y) {};
                \node[vertex,right of=y,label=above:\(c_1\)] (c1) {};
                \node[vertex,right of=c1,label=above:\(c_2\)] (c2) {};
                \node[above of=a1, label = above:\(f\)] (f){};
                \node[above of=b, label = above:\(e\)] (e){};
                \node[above of=c1, label = above:\(g\)] (g){};
    
                \begin{pgfonlayer}{background}
                    \draw[thickedge, color = yellow] (x) -- (b) -- (y);
                    \draw[mededge, color = red] (a2) -- (a1) -- (x);
                    \draw[edge, color = green] (y) -- (c1) -- (c2);
                \end{pgfonlayer}
            \end{tikzpicture}
        }
        \caption{$|f \cap g| = 0$}\label{3U:T1}
    \end{subfigure} \hfil
    
\vspace{0.5cm}

    \begin{subfigure}[h]{0.5\textwidth}
        \centering
        \scalebox{0.55}
        {
            \begin{tikzpicture}
                [every label/.append style={font=\LARGE}]
                \node[vertex,label=above:\(a_2\)] (a2) {};
                \node[vertex,right of=a2,label=above:\(a_1\)] (a1) {};
                \node[vertex,right of=a1,label=above:\(x\)] (x) {};
                \node[vertex,right of=x,label=above:\(b_1\)] (b1) {};
                \node[vertex,right of=b1,label=above:\(b_2\)] (b2) {};
                \node[vertex,below of=x,label=above:\(c_1\)] (c1) {};
                \node[vertex,below of=c1,label=above:\(c_2\)] (c2) {};
    
                \node[above of=a1, label = above:\(f\)] (f){};
                \node[above of=b1, label = above:\(e\)] (e){};
                \node[right of=c1, label = right:\(g\)] (g){};
    
                \begin{pgfonlayer}{background}
                    \draw[thickedge, color = yellow] (x) -- (b1) -- (b2);
                    \draw[mededge, color = red] (a2) -- (a1) -- (x);
                    \draw[edge, color = green] (x) -- (c1) -- (c2);
                \end{pgfonlayer}
            \end{tikzpicture}
        }
        \caption{$|f \cap g| = 1$, $f \cap g = e \cap f = e \cap g$}\label{3U:T2}
    \end{subfigure}
    \begin{subfigure}[h]{0.5\textwidth}
        \centering
        \scalebox{0.55}
        {
            \begin{tikzpicture}
                [every label/.append style={font=\LARGE}]
                \node[vertex,label=above:\(x\)] (x) {};
                \node[vertex,below left of=x,label=above:\(a\)] (a) {};
                \node[vertex,below left of=a,label=above:\(y\)] (y) {};
                \node[vertex,below right of=x,label=above:\(b\)] (b) {};
                \node[vertex,below right of=b,label=above:\(z\)] (z) {};
                \path (y) -- (z) node[vertex, midway, label=above:\(c\)] (c) {};
    
                \node[above left of=a, label = above:\(f\)] (f){};
                \node[above right of=b, label = above:\(e\)] (e){};
                \node[below of=c, label = below:\(g\)] (g){};
    
                \begin{pgfonlayer}{background}
                    \draw[thickedge, color = yellow] (x) -- (b) -- (z);
                    \draw[mededge, color = red] (x) -- (a) -- (y);
                    \draw[edge, color = green] (y) -- (c) -- (z);
                \end{pgfonlayer}
            \end{tikzpicture}
        }
        \caption{$|f \cap g| = 1$, Intersections between $e, f, g$ disjoint}\label{3U:T3}
    \end{subfigure}
\vspace{0.5cm}

    \caption{2-Matching Types when $\nu^{(k-1)} = 3$}\label{fig:3UM}
\end{figure}

    Suppose there is a matching of type~\ref{3U:T1}. If there is no edge containing $\{a_1, a_2\}$, then we are done by Lemma~\ref{lem:T_iLB}. Similarly, if there is no edge containing $\{c_1, c_2\}$, we are done. So, suppose there are some edges $f_1, g_1$ with $f_1 = \{a_1, a_2, u\}$, $g_1 = \{c_1, c_2, v\}$. If $u \not \in (e \cup g) - x$, then $e, f_1, g$ is a disconnected matching and we are done. Similarly, if $v \not \in (e \cup f) - y$, then $e, f, g_1$ is a disconnected matching and we are done. So, we may assume $u \in (e \cup g) - x$ and $v \in(e \cup f) - y$.

    If $\tu{k-1}{S_f} = 1$ or $\tu{k-1}{S_g} = 1$, we are done by Lemma~\ref{lem:pendantEdge}. Therefore, we may assume that $|S_f| > 2$ and $|S_g| > 2$. Let $f_2 \in S_f - f_1 - f$ and $g_2 \in S_g - g_1 - g$. So, $f_2 = \{a, x, u\}$, $g_2 = \{c, y, v\}$, where $a \in \{a_1, a_2\}, c \in \{c_1, c_2\}$. Since $f_2 \in S_f$ and $u \in (e \cup g) - x$, then $u$ must be in $g - e$ since otherwise, $|f_2 \cap e| = |f_2 \cap f| = 2$, a contradiction to $f_2 \in S_f$. Similarly, $v \in f - e$. Now, we obtain a 2-cover of $H$ of size exactly 6 as witnessed by $\mathcal{C} = \{\binom{e}{2}, \{u , v\}, \{f - v\}, \{g - u\}\}$. \\

    Observe that for the other cases, if there are $2$ disjoint edges in $H$, we are done. This is because either the union of their $2$-sets are a cover of $H$ or we may extend the matching to a matching of the first type or a disconnected matching.

    Next, suppose there is a matching of type~\ref{3U:T2}. By Lemma~\ref{lem:T_iLB}, $\{c_1, c_2\}$ must be contained in some edge other than $g$, say $g_1$. But then, either $g_1$ is disjoint from $e$ or $g_1$ is disjoint from $f$. In either case, we are done. \\

    In the final case, because $H$ is assumed to have no disjoint edges, it can be checked that 
    \begin{equation*}
        \mathcal{C} = \{\{x, y\}, \{x, z\}, \{y, z\}, \{x, c\}, \{y, b\}, \{z, a\}\}.
    \end{equation*}
    
    is a $2$-cover of $H$. This concludes the proof for $3$-uniform hypergraphs. \\
    
     Next, suppose $k \geq 4$ and let $H$ be a $k$-uniform hypergraph with $\vu{k-1}{H} = 3$. Let $M = \{e, f, g\}$ be a maximum $(k - 1)$-matching in $H$. Without loss of generality, suppose $|e \cap f| = k - 2$. By Lemma~\ref{lem:connected}, if $|g \cap e| \leq k - 3$ and $|g \cap f| \leq k - 3$, we are done. So, again, without loss of generality, suppose $|g \cap e| = k - 2$. We now define some notation that will be used throughout the proof. Let $S = e \cap f$, where $|S| = k - 2$ and $S' = e \cap f \cap g$. Let $e - f = \{u_1, u_2\}$, $f - e = \{v_1, v_2\}$, and $T = V(g) - e - f$. Now, $M$ will look like one of the matchings from Figure~\ref{fig:kUM}.



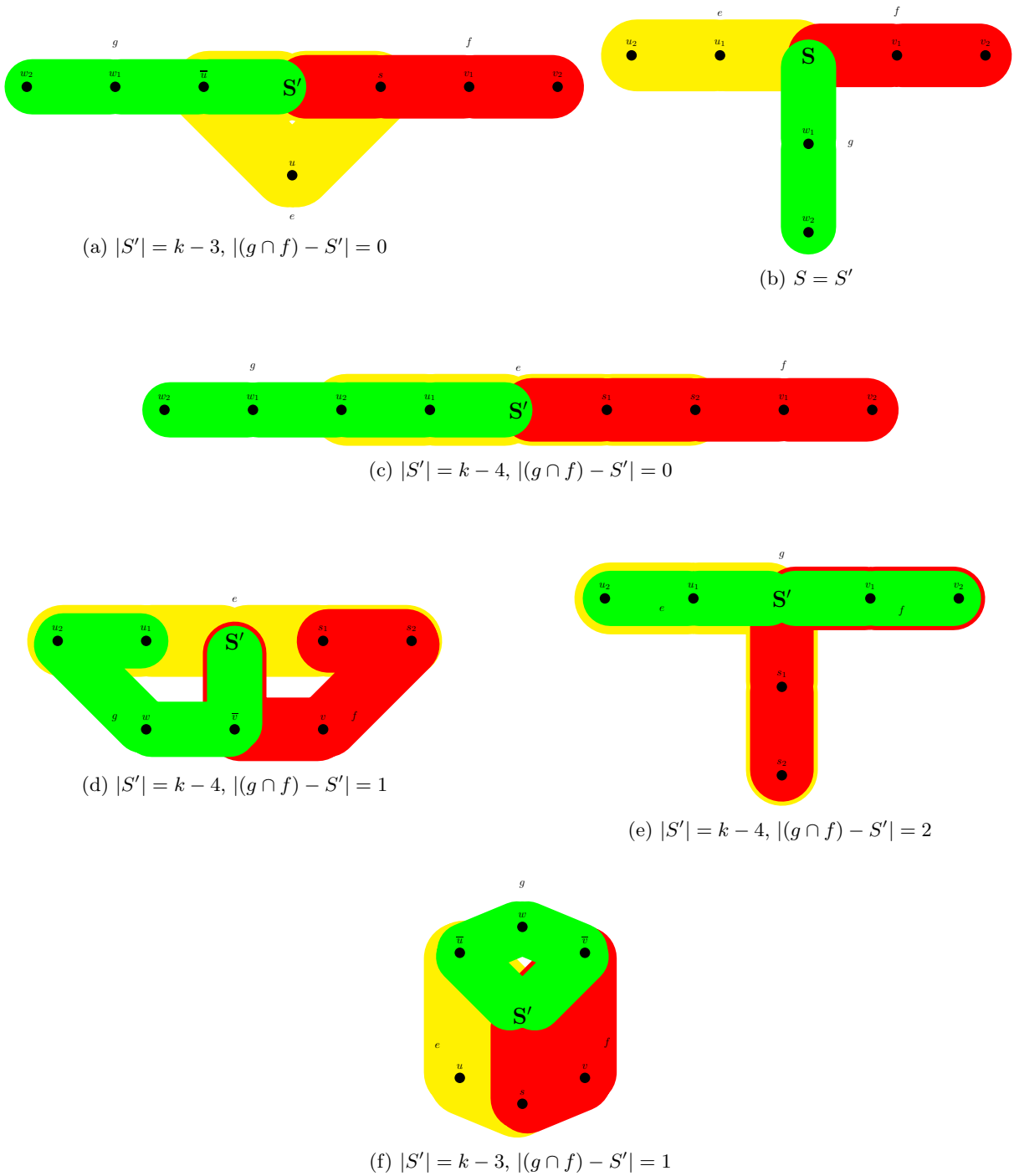
\begin{figure}
    \begin{subfigure}[h]{0.45\textwidth}
        \centering
        \scalebox{0.5}
            {
            \begin{tikzpicture}
                [every label/.append style={font=\LARGE}]
                \node[black] (Sp) {\huge $\mathbf{S'}$};
                \node[vertex, left of=Sp, label=above:\(\overline{u}\)] (ou) {};
                \node[vertex,left of=ou,label=above:\(w_1\)] (w1) {};
                \node[vertex,left of=w1,label=above:\(w_2\)] (w2) {};
                \node[vertex,right of=Sp,label=above:\(s\)] (s) {};
                \node[vertex,right of=s,label=above:\(v_1\)] (v1) {};
                \node[vertex,right of=v1,label=above:\(v_2\)] (v2) {};
                \node[vertex,below of=Sp,label=above:\(u\)] (u) {};
    
                \node[below of=u, label = below:\(e\)] (e){};
                \node[above of=v1, label = above:\(f\)] (f){};
                \node[above of=w1, label = above:\(g\)] (g){};
    
                \begin{pgfonlayer}{background}
                    \draw[thickedge, color = yellow] (Sp) -- (ou) -- (u) -- (s) -- (Sp);
                    \draw[mededge, color = red] (Sp) -- (s) -- (v1) -- (v2);
                    \draw[edge, color = green] (Sp) -- (ou) -- (w1) -- (w2);
                \end{pgfonlayer}
            \end{tikzpicture}
            }
            \caption{$|S'| = k - 3$, $|(g \cap f) - S'| = 0$}\label{kU:kMinus3T1}
    \end{subfigure}
    \hspace{1.5cm} 
    \begin{subfigure}[h]{0.45\textwidth}
        \centering
        \scalebox{0.5}
            {
            \begin{tikzpicture}
                [every label/.append style={font=\LARGE}]
                \node[black] (S) {\huge $\mathbf{S}$};
                \node[vertex,left of=S,label=above:\(u_1\)] (u1) {};
                \node[vertex, left of = u1, label=above:\(u_2\)] (u2) {};
                \node[vertex,right of=S,label=above:\(v_1\)] (v1) {};
                \node[vertex,right of=v1,label=above:\(v_2\)] (v2) {};
                \node[vertex,below of=S,label=above:\(w_1\)] (w1) {};
                \node[vertex,below of=w1,label=above:\(w_2\)] (w2) {};
    
                \node[above of=u1, label = above:\(e\)] (e){};
                \node[above of=v1, label = above:\(f\)] (f){};
                \node[right of=w1, label = right:\(g\)] (g){};
    
                \begin{pgfonlayer}{background}
                    \draw[thickedge, color = yellow] (u2) -- (u1) -- (S);
                    \draw[mededge, color = red] (S) -- (v1) -- (v2);
                    \draw[edge, color = green] (S) -- (w1) -- (w2);
                \end{pgfonlayer}
            \end{tikzpicture}
            }
            \caption{$S = S'$}\label{kU:kMinus2}
    \end{subfigure}

    \vspace{1cm}

    \begin{subfigure}[h]{\textwidth}
        \centering
        \scalebox{0.5}
            {
            \begin{tikzpicture}
                [every label/.append style={font=\LARGE}]
                \node[black] (Sp) {\huge $\mathbf{S'}$};
                \node[vertex, left of=Sp, label=above:\(u_1\)] (u1) {};
                \node[vertex, left of=u1, label=above:\(u_2\)] (u2) {};
                \node[vertex, left of=u2, label=above:\(w_1\)] (w1) {};
                \node[vertex, left of=w1, label=above:\(w_2\)] (w2) {};
                \node[vertex, right of=Sp, label=above:\(s_1\)] (s1) {};
                \node[vertex, right of=s1, label=above:\(s_2\)] (s2) {};
                \node[vertex, right of=s2, label=above:\(v_1\)] (v1) {};
                \node[vertex, right of=v1, label=above:\(v_2\)] (v2) {};
    
                \node[above of=Sp, label = above:\(e\)] (e){};
                \node[above of=v1, label = above:\(f\)] (f){};
                \node[above of=w1, label = above:\(g\)] (g){};
    
                \begin{pgfonlayer}{background}
                    \draw[thickedge, color = yellow] (u2) -- (u1) -- (Sp) -- (s1) -- (s2);
                    \draw[mededge, color = red] (Sp) -- (s1) -- (s2) -- (v1) -- (v2);
                    \draw[edge, color = green] (w2) -- (w1) -- (u2) -- (u1) -- (Sp);
                \end{pgfonlayer}
            \end{tikzpicture}
            }
            \caption{$|S'| = k - 4$, $|(g \cap f) - S'| = 0$}\label{kU:kMinus4T1}
    \end{subfigure}

    \vspace{1cm}

    \begin{subfigure}[h]{0.45\textwidth}
        \centering
        \scalebox{0.5}
            {
            \begin{tikzpicture}
                [every label/.append style={font=\LARGE}]
                \node[black] (Sp) {\huge $\mathbf{S'}$};
                \node[vertex, left of=Sp, label=above:\(u_1\)] (u1) {};
                \node[vertex,left of=u1,label=above:\(u_2\)] (u2) {};
                \node[vertex,right of=Sp,label=above:\(s_1\)] (s1) {};
                \node[vertex,right of=s1,label=above:\(s_2\)] (s2) {};
                \node[vertex,below of=Sp,label=above:\(\overline{v}\)] (ov) {};
                \node[vertex,right of=ov,label=above:\(v\)] (v) {};
                \node[vertex,left of=ov,label=above:\(w\)] (w) {};
    
                \node[above of=Sp, label = above:\(e\)] (e){};
                \node[right of=v, label = above:\(f\)] (f){};
                \node[left of=w, label = above:\(g\)] (g){};
    
                \begin{pgfonlayer}{background}
                \draw[thickedge, color = yellow] (u2) -- (u1) -- (Sp) -- (s1) -- (s2);
                    \draw[mededge, color = red] (Sp) -- (ov) -- (v) -- (s2) -- (s1);
                    \draw[edge, color = green] (Sp) -- (ov) -- (w) -- (u2) -- (u1);
                \end{pgfonlayer}
            \end{tikzpicture}
            }
            \caption{$|S'| = k - 4$, $|(g \cap f) - S'| = 1$}\label{kU:kMinus4T2}
    \end{subfigure}
    \hspace{1cm} 
    \begin{subfigure}[h]{0.45\textwidth}
        \centering
        \scalebox{0.5}
            {
            \begin{tikzpicture}
                [every label/.append style={font=\LARGE}]
                \node[black] (Sp) {\huge $\mathbf{S'}$};
                \node[vertex, left of=Sp, label=above:\(u_1\)] (u1) {};
                \node[vertex,left of=u1,label=above:\(u_2\)] (u2) {};
                \node[vertex,right of=Sp,label=above:\(v_1\)] (v1) {};
                \node[vertex,right of=v1,label=above:\(v_2\)] (v2) {};
                \node[vertex,below of=Sp,label=above:\(s_1\)] (s1) {};
                \node[vertex,below of=s1,label=above:\(s_2\)] (s2) {};
    
                \node[above of=Sp, label = above:\(g\)] (g){};
                \node[left of=u1, label = below:\(e\)] (e){};
                \node[right of=v1, label = below:\(f\)] (f){};
    
                \begin{pgfonlayer}{background}
                    \draw[thickedge, color = yellow] (u2) -- (u1) -- (Sp) -- (s1) -- (s2);
                    \draw[mededge, color = red] (s2) -- (s1) -- (Sp) -- (v1) -- (v2);
                    \draw[edge, color = green] (u2) -- (u1) -- (Sp) -- (v1) -- (v2);
                \end{pgfonlayer}
            \end{tikzpicture}
            }
            \caption{$|S'| = k - 4$, $|(g \cap f) - S'| = 2$}\label{kU:kMinus4T3}
    \end{subfigure}
    
\vspace{0.5cm}

    \begin{subfigure}[h]{\textwidth}
        \centering
        \scalebox{0.5}
            {
            \begin{tikzpicture}
                [every label/.append style={font=\LARGE}]
                \node[black] (Sp) {\huge $\mathbf{S'}$};
                \node[vertex, above of=Sp, label=above:\(w\)] (w) {};
                \node[vertex,above left of=Sp,label=above:\(\overline{u}\)] (ou) {};
                \node[vertex,above right of=Sp,label=above:\(\overline{v}\)] (ov) {};
                \node[vertex,below of=Sp,label=above:\(s\)] (s) {};
                \node[vertex,below right of=Sp,label=above:\(v\)] (v) {};
                \node[vertex,below left of=Sp,label=above:\(u\)] (u) {};
    
                \node[above left of=u, label = above:\(e\)] (e){};
                \node[above right of=v, label = above:\(f\)] (f){};
                \node[above of=w, label = above:\(g\)] (g){};
    
                \begin{pgfonlayer}{background}
                    \draw[thickedge, color = yellow] (Sp) -- (s) -- (u) -- (ou) -- (Sp);
                    \draw[mededge, color = red] (Sp) -- (ov) -- (v) -- (s) -- (Sp);
                    \draw[edge, color = green] (Sp) -- (ou) -- (w) -- (ov) -- (Sp);
                \end{pgfonlayer}
            \end{tikzpicture}
            }
            \caption{$|S'| = k - 3$, $|(g \cap f) - S'| = 1$}\label{kU:kMinus3T2}
    \end{subfigure}
    
    \vspace{0.5cm}
    
    \caption{($k - 1$)-Matching Types when $\nu^{(k-1)} = 3$, $k \geq 4$}\label{fig:kUM}
\end{figure}

    In their respective pictures, $s, s_1, s_2 \in S - S'$, $w, w_1, w_2 \in T$, $\{u, \overline{u}\} = \{u_1, u_2\}$, and $\{v, \overline{v}\} = \{v_1, v_2\}$. Throughout the proof, we will often use the result from Theorem~\ref{thm:main1} and arguments similar to the proof of the 3-uniform case. \\
    
    If we have a type~\ref{kU:kMinus3T1} matching, then observe that no edge $e' \in S_e - e$ may contain the set $\{\overline{u}, u, s\}$ since then, $e', g, f$ is a disconnected matching and we are done. Therefore, we may $(k-1)$-cover $T_e$ with three sets, namely $S' \cup A$ for each $A \in \binom{\{\overline{u}, u, s\}}{2}$. After covering $T_e$, $\nu^{(k-1)}(H - T_e) = 2$ with $M - e$ being a maximum $(k - 1)$-matching. Now, by Theorem~\ref{thm:main1}, we may find a $(k-1)$-cover of $H$ of size at most
    \begin{equation*}
        3 + 2 \left \lceil \frac{k + 1}{2} \right \rceil \leq 3 \left \lceil \frac{k + 1}{2} \right \rceil.
    \end{equation*}

    Suppose there is no type~\ref{kU:kMinus3T1} matching. If there is a type~\ref{kU:kMinus2} matching, then for all $h \in M$, there is no $h' \in S_h - h$ such that $h'$ contains $h_1, h_2$. (This is because if such an $h'$ existed, then $M - h + h'$ is a disconnected matching or a type~\ref{kU:kMinus3T1} matching.) Therefore, for each $h \in M$, we may $(k-1)$-cover $T_h$ with the sets $S + h_1$ and $S + h_2$, giving us a $(k-1)$-cover of $H$ of size 6, which is less than $3 \lceil \frac{k + 1}{2} \rceil$. \\
    
    Next, suppose there is no type~\ref{kU:kMinus3T1} or~\ref{kU:kMinus2} matching. If there is a type~\ref{kU:kMinus4T1} matching, then notice that no edge in $S_e$ contains the set $\{u_2, u_1, s_1, s_2\}$ since otherwise, we would be able to find a disconnected matching. Therefore, we may $(k-1)$-cover $T_e$ with four sets, namely $S' \cup A$ for each $A \in \binom{\{u_2, u_1, s_1, s_2\}}{3}$. If $\tu{k-1}{S_g} = 2$ or $\tu{k-1}{S_f} = 2$, we are done. Otherwise, for $h \in \{g, f\}$, we know that there exists some $h' \in S_h - h$ such that $h - e \subseteq h'$. This tells us that for all $g' \in S_g - g$, the unique vertex outside of $g' - g$ described in Lemma~\ref{lem:uniqVtx} must be $s$, where $s \in \{s_1, s_2\}$ (if not, then for any $g' \in S_g - g$ with $g - e \subseteq g'$, $M - g + g'$ is a disconnected matching). Similarly, for all $f' \in S_f - f$, the unique vertex outside of $f' - f$ described in Lemma~\ref{lem:uniqVtx} must be $u$, where $u \in \{u_1, u_2\}$. Therefore, every uncovered edge in $S_g - g$ has the form $S'' \cup \{w_1, w_2, s\}$, where $S'' \in \binom{S' \cup \{u_1, u_2\}}{k - 3}$. By Lemma~\ref{lem:trueTau}, we may cover these edges as well as $g$ with at most $\left \lceil \frac{k - 1}{2} \right \rceil$ $(k - 1)$-sets. A symmetric argument shows that we may cover the remaining uncovered edges of $S_f$ (including $f$) with at most $\left \lceil \frac{k - 1}{2} \right \rceil$ $(k - 1)$-sets. Now, we have found a cover of $H$ of size at most
    \begin{equation*}
        4 + 2 \left \lceil \frac{k - 1}{2} \right \rceil \leq 3 \left \lceil \frac{k + 1}{2} \right \rceil.
    \end{equation*}

    Now, suppose there is no type~\ref{kU:kMinus3T1} -~\ref{kU:kMinus4T1} matching and suppose there is a type~\ref{kU:kMinus4T2} matching. If $k = 4$, we cover $H$ as follows. First, we add $\{u_1, u_2, w\}, \{u_1, u_2, \overline{v}\}, \{s_1, s_2, v\}, \{s_1, s_2, \overline{v}\}$ to the cover, $\mathcal{C}$. If $\tu{k-1}{S_e} = 1$, we are done. Otherwise, there is a unique vertex $x$ outside of $e$ as described in Lemma~\ref{lem:uniqVtx} such that for all $e' \in S_e - e$, $e' - e = x$. If $x \not \in g \cup f$, then for any $e' \in S_e - e$, $M - e + e'$ is a disconnected matching. Otherwise, suppose $x \in g \cup f$ and without loss of generality, suppose $x \in g$. Then, there are at most three edges in $S_e$ that are not already covered. Namely, the edges $\{s_1, s_2, u_1, x\}, \{s_1, s_2, u_2, x\}$, and $e$. By Lemma~\ref{lem:trueTau}, we may cover these edges with two additional sets. Now, we wish to show that the edges remaining uncovered in $S_g \cup S_f$ may be covered by at most three $3$-sets. By Lemma~\ref{lem:trueTau}, either we may cover the remaining elements of $S_g$ with one $3$-set or we need to cover two edges with a unique vertex outside of $g$, which may be covered by $\left \lceil \frac{2}{2} \right \rceil = 1$ set and similarly for $S_f$. In either case, we are done. \\

    Now, suppose $k \geq 5$. We begin by adding to our cover the two $(k-1)$-sets contained in $g$ which contain $S' \cup \{u_1, u_2 \}$ and the two $(k-1)$-sets contained in $f$ which contain $S' \cup \{s_1, s_2\}$. First, we aim to cover $S_e$. Either $\tu{k-1}{S_e} = 1$ or there is a unique vertex $x$ outside of $e$ such that for all $e' \in S_e - e$, $e' - e = x$. If $x \neq \overline{v}$, then for any $e' \in S_e - e$, $M - e + e'$ is a disconnected matching. So, we may assume $x = \overline{v}$. Now, any edge $e' \in S_e - e$ which contains all of $S'$ has already been covered. Therefore, all remaining uncovered edges of $S_e - e$ have the form $S'' \cup \{u_1, u_2, v_1, v_2, \overline{v}\}$ for some $S'' \in \binom{S'}{k - 5}$. Since $e$ also remains uncovered, we are left to cover at most $k - 3$ additional edges, which by Lemma~\ref{lem:trueTau}, may be done using at most $\left \lceil \frac{k - 3}{2} \right \rceil$ $(k - 1)$-sets.

    We will now make an argument for $S_g$, which will hold true for $S_f$ by symmetry. The remaining edges of $S_g$ needing to be covered must use both $w$ and $\overline{v}$. Suppose the remaining edges of $S_g$ may not be covered by a single $(k - 1)$-set. Then, by Lemma~\ref{lem:uniqVtx}, there is a unique vertex $y$ outside of $g$ such that for all $g' \in S_g - g$, $g' - g = y$. This tells us that all edges uncovered in $S_g$ have the form $S'' \cup \{w, \overline{v}, y\}$, where $S'' \in \binom{S' \cup \{u_1, u_2\}}{k - 3}$. Specifically, there are at most $k - 2$ remaining edges to cover in $S_g$. By Lemma~\ref{lem:trueTau}, we may cover these edges with at most $\left \lceil \frac{k - 2}{2} \right \rceil$ $(k - 1)$-sets. We may make the same argument for the uncovered edges of $S_f$. All together, we have found a cover for $H$ of size:
    \begin{equation*}
        4 + \left \lceil \frac{k - 3}{2} \right \rceil + 2 \left \lceil \frac{k - 2}{2} \right \rceil \leq 3 \left \lceil \frac{k + 1}{2} \right \rceil.
    \end{equation*}
    
    Next, suppose there is no type~\ref{kU:kMinus3T1} -~\ref{kU:kMinus4T2} matching and suppose there is a type~\ref{kU:kMinus4T3} matching. We will make an argument for $S_g$, which will hold true for $S_e, S_f$ by symmetry. Suppose $\tu{k-1}{S_g} \neq 1$. Then, by Lemma~\ref{lem:uniqVtx}, there is a unique vertex $x$ outside of $g$ such that for all $g' \in S_g - g$, $g' - g = x$. Suppose $x \not \in (e \cup f) - g$. Then, for any $g' \in S_g - g$, $M - g + g'$ is either a disconnected matching or a type~\ref{kU:kMinus4T2} matching. Therefore, $x \in \{s_1, s_2\}$. Now, if any edge of $S_g - g$ contains $S'$, then this edge is actually an element of $T_g - S_g$. Therefore, every edge in $S_g - g$ has the form $S'' \cup \{u_1, u_2, v_1, v_2, x\}$, where $S'' \in \binom{S'}{k - 5}$. Since we wish to also cover $g$, there are at most $k - 3$ edges needed to be covered in $S_g$. By Lemma~\ref{lem:trueTau}, this may be done using at most $\left \lceil \frac{k - 3}{2} \right \rceil$ $(k - 1)$-sets. Similarly, the edges of $S_e$ and $S_f$ may be covered with at most $\left \lceil \frac{k - 3}{2} \right \rceil$ $(k - 1)$-sets. We are left to cover the edges which intersect more than one of $e, f, g$ in $k - 1$ vertices. We cover the edges intersecting both $g$ and $e$ in $k - 1$ vertices with the two $(k-1)$-sets contained in $e$ which contain $S' \cup \{u_1, u_2\}$. We cover the edges intersecting both $g$ and $f$ in $k-1$ vertices with the two $(k-1)$-sets contained in $f$ which contain $S' \cup \{v_1, v_2\}$. Finally, we cover the edges intersecting $e$ and $f$ with the two $(k-1)$-sets contained in $f$ which contain $S' \cup \{s_1, s_2\}$. All together, we have found a cover of $H$ of size at most
    \begin{equation*}
        3 \left \lceil \frac{k - 3}{2} \right \rceil + 6 \leq 3 \left \lceil \frac{k + 1}{2} \right \rceil.
    \end{equation*}
    
    Finally, suppose there is only a matching of type~\ref{kU:kMinus3T2}. We first show that in this case, there are no two edges with intersection size $k - 3$. For sake of contradiction, suppose there exists $h, h' \in H$ such that $|h \cap h'| = k - 3$. Let us set $A = h \cap h'$. Then, either $h, h'$ may be extended to a matching of size 3 or $h, h'$ is a maximal matching. In the first case, the extended matching must be disconnected or a matching of type~\ref{kU:kMinus3T1} or~\ref{kU:kMinus4T2}. Suppose then that $h, h'$ is a maximal matching. That is, every edge of $H$ intersects $h$ or $h'$ in $k - 1$ vertices. Because $|h \cap h'| = k - 3$, then no edge of $H$ can intersect both $h$ and $h'$ in $k - 1$ vertices. Now, we construct a suitable cover in this case. First, we cover all edges containing $A$ with the three $(k-1)$-sets contained in $h$ which contain $A$ and the three $(k-1)$-sets contained in $h'$ which contain $A$. Observe that we have also covered $h$ and $h'$. \\
    Next, let $H_h$ be the set of uncovered edges intersecting $h$ in $k - 1$ vertices and define $H_{h'}$ similarly. We will make an argument for $H_h$, which will hold true by symmetry for $H_{h'}$. First, observe that $\vu{k-1}{H_h} = 1$. Indeed, otherwise, we may find a disconnected matching of size $3$ in $H$. Also, it is the case that $\vu{k-1}{H_h \cup h} = 1$. This is because by the way $H_h$ is defined, any matching of size two in $H_h \cup h$ does not contain $h$. Now, suppose $\tu{k-1}{H_h \cup h} > 1$. Then, by Lemma~\ref{lem:uniqVtx}, there is a unique vertex $v$ outside of $h$ such that $v \in e$ for all $e \in H_h$. This means that every edge of $H_h$ has the form $(A' \cup h-A) + v$, where $A' \in \binom{A}{k - 4}$. This shows that $|H_h| \leq k - 3$ and so, by Lemma~\ref{lem:trueTau}, we may find a cover of $H_h$ of size at most $\left \lceil \frac{k - 3}{2} \right \rceil$. Similarly, $\tu{k-1}{H_{h'}} \leq \left \lceil \frac{k - 3}{2} \right \rceil$. Putting this together, we have found a cover of $H$ of size at most
    \begin{equation*}
        6 + 2 \left \lceil \frac{k - 3}{2} \right \rceil \leq 3 \left \lceil \frac{k + 1}{2} \right \rceil.
    \end{equation*}
    
    For the remainder of the proof, we may assume that no two edges intersect in exactly $k - 3$ vertices. Now, we proceed assuming that there is only a matching of type~\ref{kU:kMinus3T2}. We may cover $(T_e \cup T_f \cup T_g) - (S_e \cup S_f \cup S_g)$ with the three $(k-1)$-sets containing $S'$ and exactly two elements from $\{\overline{u}, \overline{v}, s\}$. \\
    Next, we make an argument for the uncovered edges of $S_g$, which holds true for $S_e$, $S_f$ by symmetry. Suppose $\tu{k-1}{S_g} > 1$. Then, by Lemma~\ref{lem:uniqVtx}, there is a unique vertex $x$ outside of $g$ such that for all $g' \in S_g - g$, $g' - g = x$. If $x \not \in (e \cup f) - g$, then there is an uncovered $g' \in S_g - g$ such that $M - g + g'$ is either a disconnected matching or a matching of type~\ref{kU:kMinus3T1}. Therefore, we may assume $x \in (e \cup f) - g$. This tells us all uncovered edges of $S_g$ contain $x$ and $w$. For any choice of $x$, there are at most $k - 2$ uncovered edges of $S_g$. By Lemma~\ref{lem:trueTau}, these uncovered edges of $S_e$ may be covered by at most $\left \lceil \frac{k - 2}{2} \right \rceil$ $(k - 1)$-sets. Since a symmetric argument is true for $S_e$ and $S_f$, we have found a cover of $H$ of size at most
    \begin{equation*}
        3 + 3 \left \lceil \frac{k - 2}{2} \right \rceil \leq 3 \left \lceil \frac{k + 1}{2} \right \rceil.
    \end{equation*}
\end{proof}

\section{Bounds on $g_1(k, m)$}

We begin this section with a useful definition and observation.

\begin{definition}
    Let $H$ be a $k$-uniform hypergraph and let $e \in E(H)$. For $2 \leq m \leq k - 1$, we call an $m$-set $a$ of $e$ \emph{dispensable} if for every $f \in E(H)$, $f$ intersects $e$ in some $m$-set other than $a$. Otherwise, we call $a$ \emph{indispensable}. \\

    For an indispensable $m$-set $a$ of $e$, we call any edge $f \in E(H)$ such that $f \cap e = a$ a \emph{witness} to the indispensability of $a$.
\end{definition}

\begin{observation}\label{obs:uniqWitness}
    Let $H$ be a $k$-uniform hypergraph with $m$-matching number $1$, where $\frac{k}{2} \leq m \leq k - 2$. Let $e \in E(H)$. If there is a pair of indispensable $m$-sets $a, b$ of $e$ such that $|a \cap b| = 2m - k$, there exist unique witnesses $f, g$ to $a, b$, respectively. Furthermore, we can $m$-cover $f$ and $g$ with one $m$-set.
\end{observation}

\begin{lemma}\label{lem:noMinPairs}
    Let $H$ be a $k$-uniform hypergraph with $m$-matching number $1$, $m \geq 2$. Let $e \in H$ and set $m' = \max \{0, 2m - k\}$. For any set $S \subseteq \binom{e}{m}$ of $m$-sets of $e$ with $|S| > \frac{1}{2} \binom{k}{m}$, there exists a pair $a, b \in S$ such that $|a \cap b| = m'$.
\end{lemma}

\begin{proof}
    Let $G_e$ be a graph with vertex set $\binom{e}{m}$. For $u, v \in V(G_e)$, $uv \in E(G_e)$ if and only if $|u \cap v| = m'$. Then, $G_e$ is an $\ell$-regular graph, where $\ell = \binom{k - m}{m}$ when $m' = 0$ and $\ell = \binom{m}{2m - k}$ when $m' > 0$. Observe that an independent set $I$ in $G_e$ corresponds to a set $S$ of $m$-sets of $e$ such that for any pair $a, b \in I$, $|a \cap b| \neq m'$. Using the fact that for any graph $G'$, $\alpha(G') \leq \frac{|E(G')|}{\Delta(G')}$, we have:
    \begin{equation*}
        \alpha(G_e) \leq \frac{|E(G_e)|}{\Delta(G_e)} = \frac{\left( \frac{|V(G_e)|\ell}{2} \right)}{\ell} = \frac{|V(G_e)|}{2} = \frac{1}{2} \binom{k}{m}
    \end{equation*}
    
    The result follows.
\end{proof}

We will also need the following inequality in order to prove Theorem~\ref{thm:mGEQk2}:
\begin{lemma}\label{lem:NecIneq}
    For all $k \geq 6$, $\frac{k}{2} \leq m \leq k - 2$, $0 \leq m' < m$,
    \begin{equation*}
        \binom{k}{m} > 4m - 2m' - 4
    \end{equation*}

    In particular,
    \begin{equation*}
        \binom{k}{m} - m' - 2(m - m' - 1) > \frac{1}{2} \binom{k}{m}
    \end{equation*}
\end{lemma}

\begin{proof}
    Fix $k \geq 6$, $\frac{k}{2} \leq m \leq k - 2$, and $0 \leq m' < m$. First, observe that
    \begin{equation*}
        4m - 2m' - 4 \leq 4m - 4 \leq 4(k - 2) - 4 = 4(k - 3)
    \end{equation*}

    On the other hand, we have:
    \begin{equation*}
        \binom{k}{m} \geq \binom{k}{k - 2} = \binom{k}{2}
    \end{equation*}

    Now, it is left to show the following inequality
    \begin{equation*}
        \binom{k}{2} - 4(k - 3) = \frac{1}{2}(k^2 - 9k + 24) > 0
    \end{equation*}
    
    Let $f(k) = \frac{1}{2}(k^2 - 9k + 24)$. It can be checked that $f(6) = 3 > 0$. Furthermore, $f'(k) > 0$ for all $k \geq 5$. So, $f$ is increasing for all $k \geq 5$ and therefore, $f(k) > 0$ for all $k \geq 6$. We obtain the second part of the lemma by rearranging the inequality.
\end{proof}

With the help of the above two lemmas, we are able to prove Theorem~\ref{thm:mGEQk2}.

\begin{proof}[Proof of Theorem~\ref{thm:mGEQk2}]
    Let $k \geq 6$, $\frac{k}{2} \leq m \leq k - 2$, and let $H$ be a $k$-uniform hypergraph with $m$-matching number $1$. Fix $e \in E(H)$ with the most dispensable $m$-sets. Observe that for any non-witnessing edge $f \in E(H)$, $f$ contains at least $m + 1$ $m$-sets of $e$. If $e$ has at least $m$ dispensable $m$-sets, then we may delete any $m$ of them and obtain an $m$-cover of $H$ with the remaining $m$-sets of $e$. Suppose then that $e$ has $m' < m$ dispensable $m$-sets. Denote the set of dispensable $m$-sets of $e$ by $S$. So the number of indispensable sets is $\binom{k}{m} - m'$. We wish to find an $m$-cover of size at most $\binom{k}{m} - m = (\binom{k}{m} - m') - (m - m')$. We do this by deleting $S$ from $\binom{e}{m}$ and then finding $m - m'$ pairs of indispensable $m$-sets $a_i, b_i \in \binom{e}{m} - S$ such that $|a_i \cap b_i| = 2m - k$ for $1 \leq i \leq m - m'$. \\
    
    Note that while $0 \leq i - 1 \leq m - m' - 1$, $\binom{k}{m} - m' - 2i \geq \binom{k}{m} - m' - 2(m - m' - 1)$. Set $i = 0$ and $S' = \binom{e}{m} - S$. While $i \leq m - m' - 1$, by Lemmas~\ref{lem:noMinPairs} and ~\ref{lem:NecIneq}, there exists a pair of indispensable $m$-sets of $e$, $a_i, b_i$, with witnessing edges $f_i, g_i$, respectively, such that $|a_i \cap b_i| = 2m - k$. We may cover $f_i, g_i$ with the $m$-set $x_i = (a_i \cap b_i) \cup (f_i - e)$. Note that every non-witnessing edge other than $e$ contains either at most one of $a$ and $b$. Delete $a_i, b_i$ from $S'$, increase $i$ by 1, and repeat. Now, we have the following $m$-cover $\mathcal{C}$ of $H$:
    \begin{equation*}
        \mathcal{C} = \left( \binom{e}{m} - S' \right) \cup \left(\bigcup_{i = 0}^{m - m' - 1} \{x_i\} \right) = \left( \binom{e}{m} - \left(S \cup \left(\bigcup_{i = 0}^{m - m' - 1} \{a_i, b_i\} \right) \right) \right) \cup \left(\bigcup_{i = 0}^{m - m' - 1} \{x_i\} \right)
    \end{equation*}

    Now, we can compute $|\mathcal{C}|$:
    \begin{equation*}
        \begin{split}
            |\mathcal{C}|
                &=  \binom{k}{m} - (m' + 2(m - m')) + (m - m') \\
                &=  \binom{k}{m} - m' - 2(m - m') + (m - m') \\
                &=  \binom{k}{m} - m' - (m - m') \\
                &=  \binom{k}{m} - m
        \end{split}
    \end{equation*}

    Therefore, $g_1(k, m) \leq \binom{k}{m} - m$ for all $\frac{k}{2} \leq m \leq k - 2$.
\end{proof}

Next, we improve the previous upper bound for $g_1(5, 2)$ following a similar argument as the above proof.

\begin{proof}[Proof of Theorem~\ref{thm:g52}]
    Let $H$ be a $5$-uniform hypergraph with $2$-matching number 1. Let $r = \max \{|e \cap f| : e, f \in E(H)\}$. If $r \geq 3$, then letting $e, f \in E(H)$ such that $|e \cap f| = r$, we may cover $H$ with the $2$-sets $\binom{e \cap f}{2}$ together with the $2$-sets containing exactly one element from $e - f$ and one element from $f - e$. This gives a cover of size 7. Suppose then that $r = 2$. That is, every edge intersects every other edge in exactly two vertices. Let $e$ be an edge with the most dispensable sets. Observe that for any dispensable set $a$ of $e$, there is no edge intersecting $e$ at $a$. If $e$ has at least $3$ dispensable sets, then we are done. Otherwise, we may assume $e$ has $m' \leq 2$ dispensable sets and therefore, $10 - m' \geq 8$ indispensable sets. Denote the set of dispensable sets by $S$. Observe that for any pair $a, b$ of indispensable $2$-sets of $e$ with $|a \cap b| = 0$, there exist unique witnesses $f, g$ of $a, b$, respectively. Let $S' = \binom{e}{2} - S$. So, $|S'| = 10 - m'$. Now, by Lemma~\ref{lem:noMinPairs}, we may find at least $\left \lceil \frac{|S'| - 5}{2}\right \rceil = \left \lceil \frac{5 - m'}{2}\right \rceil$ pairs of indispensable $2$-sets, $a_i, b_i$ for $1 \leq i \leq \left \lceil \frac{5 - m'}{2}\right \rceil$ such that $|a_i \cap b_i| \geq 2$ with witnesses $f_i, g_i$, respectively. For $1 \leq i \leq \left \lceil \frac{5 - m'}{2}\right \rceil$, we may $2$-cover $f_i, g_i$ with $f_i \cap g_i$. Now, we have the following $2$-cover of $H$:
    \begin{equation*}
        \mathcal{C} = \left(S' - \bigcup_{i = 1}^{\left \lceil \frac{5 - m'}{2}\right \rceil} \{a_i, b_i\}\right) \cup \bigcup_{i = 1}^{\left \lceil \frac{5 - m'}{2}\right \rceil} (f_i \cap g_i)
    \end{equation*}

    The size of this $2$-cover is:
    \begin{equation*}
        |\mathcal{C}| = \left((10 - m') - 2 \cdot \left \lceil \frac{5 - m'}{2}\right \rceil \right) + \left \lceil \frac{5 - m'}{2}\right \rceil = 10 - m' - \left \lceil \frac{5 - m'}{2}\right \rceil \leq 7
    \end{equation*}
\end{proof}

We next improve the bound given by Theorem~\ref{thm:mGEQk2} for the case when $m = k - 2$. We will need the following lemma:
\begin{lemma}\label{lem:nearPerfMatching}
    Let $k \geq 5$ and let $G$ be a graph with vertex set $\binom{[k]}{k - 2}$ and for $A, B \in V(G)$, $AB \in E(G)$ if and only if $|A \cap B| = k - 4$. Then, $G$ has a perfect matching if $\binom{k}{2}$ is even and $G$ has a matching with one unsaturated vertex when $\binom{k}{2}$ is odd.
\end{lemma}
\begin{proof}
    By Theorem 1.2 from~\cite{ye18}, $G$ has a maximum matching such that any pair of unsaturated vertices have no common neighbors. Therefore, if every pair of vertices have a common neighbor, we are done. When $k \geq 6$, by inclusion-exclusion, it is easy to check that for any $x, y \in V(G)$, $|N(x) \cap N(y)| > 0$. When $k = 5$, let $M$ be a maximum matching of $G$ such that any pair of unsaturated vertices have no common neighbors. Suppose $A, B$ are unsaturated by $M$. Then, $AB \not \in E(G)$ as this would contradict that $M$ is a maximum matching. This implies that $|A \cap B| = 2$. But then, the vertex $C = \{A - B, B - A, [k] - (A \cup B)\}$ is a common neighbor of $A$ and $B$, a contradiction.
\end{proof}

\begin{lemma}\label{lem:perfMatching}
    Let $k \geq 5$ and let $H$ be a $k$-uniform hypergraph with $\vu{k-2}{H} = 1$. If there exists an edge that intersects every other edge in exactly $k-2$ vertices, then
    \begin{equation*}
        \tu{k-2}{H} \leq \left \lceil \frac{\binom{k}{k-2} + 1}{2} \right \rceil = \left \lceil \frac{\binom{k}{2} + 1}{2} \right \rceil
    \end{equation*}
\end{lemma}

\begin{proof}
   Let $k \geq 5$ and let $H$ be a $k$-uniform hypergraph with $\vu{k-2}{H} = 1$. Suppose there exists an edge $e \in E(H)$ that intersects every other edge in exactly $k-2$ vertices. Using the graph $G_e$ from Lemma~\ref{lem:noMinPairs} which satisfies the properties of the graph in Lemma~\ref{lem:nearPerfMatching}, there exists a matching $M$ of $G_e$ of size $\left \lfloor \frac{|V(G_e)|}{2} \right \rfloor$.

    For each $uv \in M$, if there are witnessing edges $f_u$, $f_v$ of $u$ and $v$, respectively, these witnessing edges are unique and their intersection has size exactly $k-2$. We may cover this pair of edges with the $(k-2)$-set $f_u \cap f_v$. If there is only one of the two witnessing edges, say $f_u$, then $v$ is a dispensable $(k-2)$-set and we may cover all edges intersecting $e$ in $u$ by the $(k-2)$-set $u$. Doing this for all edges of $M$, we arrive at collection of $(k-2)$-sets covering all edges of $H - e$ with the exception of the witnessing edges of at most one $(k-2)$-set. We may cover the remaining edges with at most $1$ $(k-2)$-set, giving a $(k-2)$-cover of $H$ of size
    \begin{equation*}
        |M| + 1 = \left \lfloor \frac{\binom{k}{k-2}}{2} \right \rfloor + 1 = \left \lceil \frac{\binom{k}{2} + 1}{2} \right \rceil
    \end{equation*}
\end{proof}

We are now ready to prove Theorem~\ref{thm:mEqkMinus2}.

\begin{proof}[Proof of Theorem~\ref{thm:mEqkMinus2}]
    We will prove the odd and even case separately by induction. First, suppose $k$ is odd. It is not hard to show that $g_1(3, 1) = 3 = \frac{1}{4}(k^2 + 3)$. Now, let $H$ be a $k$-uniform hypergraph with $k \geq 5$, $k$ odd, where $\vu{k-2}{H} = 1$. Furthermore, we will assume $g_1(k - 2, k - 4) \leq \frac{1}{4}((k - 2)^2 + 3)$. If there is an edge $e$ of $H$ such that every other edge of $H$ intersects $e$ in exactly $k-2$ vertices, then by Lemma~\ref{lem:perfMatching}, we may find an $(k-2)$-cover of $H$ of size $\left \lceil \frac{\binom{k}{k-2} + 1}{2} \right \rceil = \left \lceil \frac{\binom{k}{2} + 1}{2} \right \rceil \leq \frac{1}{4}(k^2 + 3)$. \\
    
    Suppose then that there is a pair of edges $e, f$ such that $|e \cap f| = k - 1$. Let us denote $e \cap f$ by $S$ and suppose $e - S = u, f - S = v$. Observe that all edges intersect $S$ in at least $k - 3$ vertices. We may $(k-2)$-cover all edges intersecting $S$ in at least $k - 2$ vertices by the $k - 1$ $(k-2)$-sets $\binom{S}{k - 2}$. Now, observe that the uncovered edges all intersect $S$ in $k - 3$ vertices. Therefore, they must contain both $u$ and $v$ since $H$ has $(k-2)$-matching number $1$. Take $H'$ to be the $(k - 2)$ uniform hypergraph with vertex set $V(H) - \{u, v\}$ and edge set $E(H') = \{g - \{u , v\} : g \in E(H), |g \cap S| = k - 3\}$. Now, $H'$ has $(k - 4)$ matching number 1. Otherwise, there exist edges $h_1', h_2' \in H'$ such that $|h_1 \cap h_2| \leq k - 5$. But then, setting $h_1 = h_1' \cup \{u, v\}, h_2 = h_2' \cup \{u, v\}$, we find that $h_1, h_2$ is a $(k - 2)$-matching in $H$, a contradiction. By induction, we have:
    \begin{equation*}
        \tau^{(k-2)}(H') \leq g_1(k - 2, k - 4) \leq \frac{1}{4}((k - 2)^2 + 3)
    \end{equation*}
    
    Letting $C'$ be a $(k - 4)$ cover of $H'$ of size $\tau^{(k-2)}(H')$, then the following is a cover of $H$:
    \begin{equation*}
        C = \{T \cup \{u, v\} : T \in C'\} \cup \binom{S}{k - 2}
    \end{equation*}

    We compute the size of $C$ to be:
    \begin{equation*}
        |C| = \tau^{(k-2)}(H') + (k - 1) \leq \frac{1}{4}((k - 2)^2 + 3) + (k - 1) = \frac{1}{4}(k^2 + 3).
    \end{equation*} \\

    The proof for $k$ even is almost the exact same. We include it here for completeness. Suppose $k$ is now even. It was shown in~\cite{az20} that $g_1(4, 2) = 4 = \frac{1}{4}4^2$. Now, let $H$ be a $k$-uniform hypergraph with $k \geq 6$, $k$ even, where $\vu{k-2}{H} = 1$. We will assume $g_1(k - 2, k - 4) \leq \frac{1}{4}(k - 2)^2$. If there is an edge $e$ of $H$ such that every other edge of $H$ intersects $e$ in exactly $k - 2$ vertices, then by Lemma~\ref{lem:perfMatching}, we may find an $(k - 2)$-cover of $H$ of size $\left \lceil \frac{\binom{k}{k - 2} + 1}{2} \right \rceil = \left \lceil \frac{\binom{k}{2} + 1}{2} \right \rceil \leq \frac{1}{4}k^2$. \\
    
    Suppose then that there is a pair of edges $e, f$ such that $|e \cap f| = k - 1$. Let us denote $e \cap f$ by $S$ and suppose $e - S = u, f - S = v$. Observe that all edges intersect $S$ in at least $k - 3$ vertices. We may $(k - 2)$-cover all edges intersecting $S$ in at least $k - 2$ vertices by the $k - 1$ $(k - 2)$-sets $\binom{S}{k - 2}$. Now, observe that the uncovered edges all intersect $S$ in $k - 3$ vertices. Therefore, they must contain both $u$ and $v$ since $H$ has $(k - 2)$-matching number $1$. Take $H'$ to be the $(k - 2)$ uniform hypergraph with vertex set $V(H) - \{u, v\}$ and edge set $E(H') = \{g - \{u , v\} : g \in E(H), |g \cap S| = k - 3\}$. Now, $H'$ has $(k - 4)$-matching number 1. Otherwise, there exist edges $h_1', h_2' \in H'$ such that $|h_1 \cap h_2| \leq k - 5$. But then, setting $h_1 = h_1' \cup \{u, v\}, h_2 = h_2' \cup \{u, v\}$, we find that $h_1, h_2$ is a $(k - 2)$-matching in $H$, a contradiction. By induction, we have:
    \begin{equation*}
        \tau^{(k-2)}(H') \leq g_1(k - 2, k - 4) \leq \frac{1}{4}(k - 2)^2
    \end{equation*}
    
    Letting $C'$ be a $(k - 4)$ cover of $H'$ of size $\tau^{(k-2)}(H')$, then the following is a cover of $H$:
    \begin{equation*}
        C = \{T \cup \{u, v\} : T \in C'\} \cup \binom{S}{k - 2}
    \end{equation*}

    We compute the size of $C$ to be:
    \begin{equation*}
        |C| = \tau^{(k-2)}(H') + (k - 1) \leq \frac{1}{4}(k - 2)^2 + (k - 1) = \frac{1}{4}k^2.
    \end{equation*}
\end{proof}

\section{Fractional Results}

We begin this section by proving Theorem~\ref{thm:fracG2kk}:

\begin{proof}[Proof of Theorem~\ref{thm:fracG2kk}]
    Let $k \geq 2$ and $H$ be a $2k$-uniform hypergraph with $k$-matching number 1 and take $e \in H$. Begin by assigning every $m$-set contained in $e$ a weight of $\frac{1}{k + 1}$. In doing this, every edge intersecting $e$ in at least $k + 1$ vertices is fractionally $k$-covered. The remaining uncovered edges intersect $e$ in exactly $k$ vertices and currently have weight $\frac{1}{k + 1}$. Observe that for any $k$-set $S$ of $e$, there is a unique $k$-set $T$ of $e$ such that $S \cup T = e$ and $S \cap T = \emptyset$. There are exactly $\frac{1}{2} \binom{2k}{k}$ such pairs of $k$-sets of $e$. Let us label them as $\{(S_i, T_i) : 1 \leq i \leq \frac{1}{2} \binom{2k}{k}\}$. Now, for each pair $S_i, T_i$, either there is a unique pair of edges $f, g$ intersecting $S_i, T_i$, respectively or there are multiple edges intersecting one of these $k$-sets and no edges intersecting the other $k$-set. In either case, we may find a single $k$-set and assign it weight $\frac{k}{k + 1}$ in order to fractionally $k$-cover all uncovered edges intersecting $e$ at $S_i$ and $T_i$. Now, we have covered all edges with a total weight of:
    \begin{equation*}
        \frac{1}{k + 1} \binom{2k}{k} + \frac{k}{k + 1} \frac{\binom{2k}{k}}{2} = \left( \frac{1}{k + 1} + \frac{k}{2(k + 1)} \right) \binom{2k}{k} = \left( \frac{1}{2} + \frac{1}{2(k+1)} \right) \binom{2k}{k}.
    \end{equation*}
\end{proof}

We may obtain bounds on $h^*(k,m)$ from $g_1^*(k,m)$ using the following lemma. This generalizes the upper bound proof strategy of Proposition 14 in~\cite{basit22} to work for all choices of $k$ and $m$.
\begin{lemma}\label{lem:hStarFromgStar}
    For all $2 \leq m < k$, we have $h^*(k, m) \leq \frac{1}{2} \left( \binom{k}{m} + g_1^*(k, m) \right)$.
\end{lemma}

\begin{proof}
    Let $H$ be a $k$-uniform hypergraph and fix $2 \leq m < k$. Suppose $H$ has $m$-matching number $\nu$ and let $M = \{e_1, \dots, e_{\nu}\}$ be a maximum $m$-matching in $H$. Begin by assigning weight $1/2$ to all of the $m$-sets in $\bigcup_{i = 1}^{\nu} \binom{e_i}{m}$. Any edge which intersects at least 2 edges of the matching in $m$ vertices is now fractionally $m$-covered as well as any edge which intersects a matching edge in more than $m$ vertices. The uncovered edges now intersect exactly 1 matching edge in exactly $m$ vertices. For $1 \leq i \leq \nu$, let $S_{e_i} = \{f \in H : |f \cap e_i| = m \text{ and $f$ is uncovered}\}$. Clearly, all uncovered edges are contained in some $S_{e_i}$. Furthermore, for any $i$, the subgraph of $H$ with edge set $S_{e_i}$ has $m$-matching number 1. Otherwise, we may find an $m$-matching of $H$ of size larger than $M$. So, for each $i$, we may fractionally $m$-cover the uncovered edges in $S_{e_i}$ with a total weight of at most $\frac{1}{2} g_1^*(k, m)$ (We only need $\frac{1}{2} g_1^*(k, m)$ since each $m$-set of a matching edge was intially given a weight of $\frac{1}{2}$). Now, we have fractionally $m$-covered $H$ with a total weight of at most $\frac{1}{2}\left( \binom{k}{m} + g_1^*(k, m) \right) \nu$, giving us
    \begin{equation*}
        h^*(k, m) \leq \frac{1}{2} \left( \binom{k}{m} + g_1^*(k, m) \right).
    \end{equation*}
\end{proof}

As mentioned in the introduction, using Lemma~\ref{lem:hStarFromgStar} together with Theorem~\ref{thm:fracG2kk}, we obtain Corollary~\ref{cor:fracH2kk}.


    
    

Lastly, we improve the upper bound on $g_1^*(k, k - 2)$ by proving Theorem~\ref{thm:fracGkkMinus2}:

\begin{proof}[Proof of Theorem~\ref{thm:fracGkkMinus2}]
    Let $H$ be a $k$-uniform hypergraph with $(k-2)$-matching number 1. If there exists some edge $e$ of $H$ such that every other edge of $H$ intersects $e$ in $k - 1$ vertices, then assigning weight $\frac{1}{k-1}$ to every $(k - 2)$-set of $e$, we obtain a fractional $(k - 2)$-cover of size $\frac{k}{2}$. Otherwise, we may find two edges $e, f$ of $H$ such that $|e \cap f| = k - 2$. Let $S = e \cap f$. Then, for any other edge $g \in H - e - f$, $|g \cap S| \in \{k - 2, k - 3, k - 4\}$. We fractionally cover all edges intersecting $S$ in $k - 2$ vertices (including $e, f$) by assigning weight $1$ to $S$. Now, the edges which intersect $S$ in $k - 3$ vertices also intersect $e - S$ and $f - S$ in at least $1$ vertex. Assigning weight $\frac{1}{k - 3}$ to every $(k - 2)$-set of the form $S' \cup \{x, y\}$, where $S' \in \binom{S}{k - 4}$, $x \in e - S$, $y \in f - S$, we fractionally $(k - 2)$-cover all edges intersecting $S$ in $k - 3$ vertices. Also, all edges intersecting $S$ in $k - 4$ vertices are partially covered (each have weight $\frac{4}{k - 3}$). Now, for every edge $g$ intersecting $S$ in $k - 4$ vertices, $(e \cup f) - S \subseteq g$. So, assigning weight $\left( 1 - \frac{4}{k - 3}\right) \frac{1}{\binom{k - 4}{2}}$ to every $(k - 2)$-set of the form $S'' \cup ((e \cup f) - S)$, where $S'' \in \binom{S}{k - 6}$, we fractionally $(k - 2)$-cover the edges intersecting $S$ in $k - 4$ vertices and we have now covered all edges of $H$. The weight of this cover is:
    \begin{equation*}
        \begin{split}
            1 + \frac{1}{k - 3} \left( 4\binom{k - 2}{k - 4} \right) + \left( 1 - \frac{4}{k - 3}\right) \frac{1}{\binom{k - 4}{2}} \binom{k - 2}{k - 6}
                &=  1 + \frac{4 \binom{k - 2}{2}}{k - 3} + \frac{k - 7}{k - 3} \frac{1}{\binom{k - 4}{2}} \binom{k - 2}{4} \\
                &=  1 + 2(k - 2) + \frac{k - 7}{6(k - 3)}\binom{k - 2}{2} \\
                &\leq \frac{1}{6} \binom{k - 2}{2} + 2k - 3.
        \end{split}
    \end{equation*}
\end{proof}

\section*{Acknowledgements}
The author would like to thank Shira Zerbib for helpful suggestions and discussions throughout the development of this paper.
\bibliography{bibfile}
\bibliographystyle{abbrv}

\end{document}